\documentstyle[12pt]{article}
\title{\bf Impulsive Stabilization of Linear Delay Differential
Equations}

\author
{L. Berezansky $^{1}$
\\Ben-Gurion University of the Negev, \\
Department of Mathematics and Computer Science, \\
Beer-Sheva 84105, Israel, \\
E. Braverman  \\
Technion - Israel Institute of Technology,\\
Department of Mathematics, Haifa 32000, Israel }

\begin{document}
\maketitle

\footnotetext[1]{Supported by
Israel Ministry of Absorption and by Israel Ministry of Science and
Technology}

\begin{abstract}
The paper is concerned with stabilization of a
scalar delay differention equation
$$
{\dot x} (t) - \sum_{k=1}^m A_k(t)x[h_k(t)] = 0,~t\geq 0,~
x(\xi)=\varphi (\xi), \xi <0,
$$
by introducing impulses in certain moments of time
$$ x(\tau_j) = B_j x(\tau_j -0), ~j=1,2, \dots ~. $$

Explicit stability results are presented both for
the equation with positive coefficients and for the equation with $A_k$
being of arbitrary sign.
\end{abstract}

\section{Introduction}

{}~~~~~The stabilization problem is among the basic ones in the
control theory.
Introducing of impulses, i.e. jumps of a solution in successive
moments of time, is one of stabilization methods.
For ordinary differential equations there exist explicit formulas
connecting solutions of impulsive and non-impulsive equations
[1,2].
So in this case one easily obtains conditions on impulses
providing that an unstable system becomes stable.
For instance the results of the monograph [1] can be applied
to constructing of such conditions.

Delay differential equations describe a wide class of
practical models.
The impulsive stabilization problem is important for these
equations as well.
However for impulsive delay differential equations there
are no explicit stability conditions applicable to stabilization.
The results obtained in [2-4] presuppose that the equation
without impulses is stable.

This paper is concerned with explicit impulsive stabilization results
for a scalar delay differential equation.
The solution representation formula for an impulsive delay
differential equation from [5] is intensively exploited.
It should be noted that solution representation formulas
are of increasing significance in the stability investigation
of functional differential equations (see, for example,
the recent monograph [6]).

The paper is organized as follows.

1. For equations with positive coefficients
(i.e. unstable without impulses)
stability conditions are presented which
connect coefficients and impulsive conditions
in each interval between impulses.

2. "Uniform stabilization" results are given
not depending on delays.

3. We demonstrate that stabilization is also possible
without definite sign of coefficients.

\section{Preliminaries}

{}~~~~~We consider a scalar linear delay differential equation
\begin{eqnarray}
\dot{x} (t) - \sum_{k=1}^m {A_k (t) x[h_k(t)]} =
r(t),~ t \geq 0,
{}~~~~x(\xi ) = \varphi (\xi),~ \xi < 0,
\end{eqnarray}
\begin{eqnarray}
x(\tau_j) = B_j x(\tau_j - 0) , ~j=1,2, \dots ,
\end{eqnarray}
under the following assumptions

(a1) $ 0 = \tau_0 < \tau_1 < \tau_2 < \dots $
are fixed points,  $\lim_{j \rightarrow \infty} \tau_j = \infty $;

(a2) $A_k, r , ~k=1, \dots ,m $ are Lebesgue measurable
essentially
bounded in any finite segment
$[0,b]$ functions,
$ ~B_j \in {\bf R},~ j=1, \dots $, ${\bf R}$ is a real axis;

(a3) $h_k: [0, \infty) \rightarrow {\bf R}$ are Lebesgue
measurable functions, $h_k(t) \leq t$ ;

(a4) $\varphi: (- \infty,0) \rightarrow {\bf R}$
is a Borel measurable bounded function.

\underline{\bf Definition.}
A function $x: [0, \infty) \rightarrow {\bf R}$ absolutely
continuous in each $[\tau_j, \tau_{j+1} ) $
is {\bf a solution} of the impulsive equation (1),(2),
if for $t\neq\tau_j$ it satisfies (1) and for $t=\tau_j$
it satisfies (2).

\underline{\bf Definition.}
For each $s \geq 0$
a solution $X(t,s)$ of the equation
\begin{eqnarray}
\dot{x} (t) - \sum_{k=1}^m {A_k (t) x[h_k(t)]} =
0,~ t \geq s,
{}~~~~x(\xi ) = 0 ,~ \xi < s,
\end{eqnarray}
\begin{eqnarray}
x(\tau_j) = B_j x(\tau_j - 0) , ~\tau_j > s,
\end{eqnarray}
satisfying $x(s) = 1$
is said to be {\bf a fundamental function }
of the equation (1),(2).
We assume $X(t,s) = 0, ~0 \leq t < s$.

Similarly a fundamental function $C(t,s)$ of the equation (1)
without impulses is defined.

\newtheorem{uess}{Lemma}
\begin{uess}
{\bf [5]} Suppose (a1)-(a4) are satisfied.
Then there exists one and only one
solution of the problem (1) , with the initial
value
$x(0) = \alpha _0
$
and impulsive conditions
$$ x(\tau_j) = B_j x(\tau_j - 0) + \alpha_j$$
and it can be presented as
\begin{eqnarray}
x(t) = X(t,0)x(0) + \int_0^t {X(t,s) r(s) ds} + \nonumber \\
+ \sum_{k=1}^m {\int_0^t {X(t,s) A_k(s) \varphi [h_k(s)] ds}}
+ \sum_{j=1}^{\infty} X(t, \tau_j) \alpha_j .
\end{eqnarray}
Here $\varphi (\zeta) = 0,$ if $\zeta \geq 0. $
\end{uess}

\underline{\bf Remark.}
Under the hypotheses of Lemma 1 for the equation (1) without
impulses the solution representation formula
\begin{eqnarray}
x(t) = C(t,0)x(0) + \int_0^t {C(t,s) r(s) ds}
+ \sum_{k=1}^m {\int_0^t {C(t,s) A_k(s) \varphi [h_k(s)] ds}} .
\end{eqnarray}
is also valid.
Here $\varphi (\xi) = 0$,  if $\xi \geq 0$.
Besides, when changing the initial point $t=0$
by an arbitrary initial point $t_0 > 0$
we obtain a the representations similar to (5) and (6),
with the same fundamental functions $X(t,s)$ and $C(t,s)$.

\underline{\bf Definition.}
The equation (1),(2) is said to be {\bf stable}
if there exists $N>0$ such that for a solution $x$
of this equation, with $f \equiv 0$, the following estimate holds
$$ \mid \! x(t) \! \mid  < N \left( \mid \! x(0) \! \mid
+ \sup_{t<0} \mid \! \varphi (t) \! \mid \right) .$$

The equation (1),(2) is said to be {\bf exponentially stable}
if there exist $N>0$ and $ \lambda > 0$ such that for
a solution $x$ of this equation, with $f \equiv 0$,
the following estimate holds
$$ \mid \! x(t) \! \mid < N e^{- \lambda t} \left(
\mid \! x(0) \! \mid + \sup_{t<0} \mid \! \varphi (t) \! \mid
\right) . $$

If $A_k (t) \geq 0, ~k=1, \dots, m,$ and for at least
one of indices $k$  $A_k (t) \geq \eta > 0$
then it is well known [7] that the non-impulsive equation (1)
is unstable.
Our objective is to derive conditions on impulses $B_j$
providing (1),(2) is stable.

The following assertion is a special case of the result
obtained in the paper [8].
However it has a simple proof that we present here.

\begin{uess}
Suppose (a2)-(a4) hold and $A_k \geq 0$.

Then the fundamental function $C(t,s)$ of the equation (1)
without impulses is positive: $C(t,s) > 0, ~0 \leq s \leq t <
\infty $.
Besides, if $\varphi (t) \geq 0, ~f(t) \geq 0 $
then the inequality $x(0) > 0$ implies $x(t) > 0, ~ t \geq 0$.
\end{uess}

{\bf Proof.}
Since $C(t,s)$ is a solution of (3) with the initial condition
$C(t,s) = 1$, then
\begin{eqnarray}
C(t,s) = 1 + \sum_{k=1}^m \int_s^t A_k (\xi) C(h_k (\xi), s)~d\xi,
t \geq s, \label{e7} \\
C(h_k (\xi), s) = 0, \mbox{~if~} h_k (\xi) < s. \nonumber
\end{eqnarray}

We prove $C(t,s) >0, ~0 \leq s \leq t < \infty$.
Assume the contrary.
Since $C(s,s) = 1$ then the continuity of $C(t,s)$ in $t$
yields $C(t,s) > 0$ in a certain interval $t \in [s, \zeta_s)$.
Let $t_s$ be the first point such that $C(t_s, s) = 0$.
Then (\ref{e7}) gives
$$ 0 = 1 + \sum_{k=1}^m \int_s^{t_s} A_k (\xi) C(h_k(\xi),s) d\xi . $$

Since the integral in the right-hand side is not negative
then the contradiction obtained proves $C(t,s)>0,
{}~0\leq s \leq t < \infty$.

The second statement of the theorem immediately follows
from the representation (6).
\vspace{2 mm}

{\bf Corollary.}
Let (a1)-(a4) hold, $A_k \geq 0,~B_j \geq 0$.

Then the fundamental function of the problem (1),(2)
is positive: $X(t,s) > 0, ~0 \leq s \leq t < \infty$.

{\bf Proof.}
Let $s>0$ be fixed and $\tau_{j_s}$ be the smallest
of $\tau_j$ such that $\tau_j >s$.
For $t \in [s, \tau_{j_s}) ~X(t,s)=C(t,s)$.
Let $t \in [\tau_{j_s},\tau_{j_s +1}]$.
Then $X(t,s)$ is a solution of the problem
$$\dot{x}(t)= \sum_{k=1}^m A_k (t) x[h_k(t)],
{}~t \in [\tau_{j_s}, \tau_{j_s+1}), $$
\begin{eqnarray}
x(\xi)=C(\xi,s),~\xi < \tau_{j_s},
{}~x(\tau_{j_s})=B_{j_s} C(\tau_{j_s},s). \label{e8}
\end{eqnarray}
By Lemma 2 with
$t=0$ being changed
by $t=t_0~~$
$x(t)=X(t,s) > 0$.
By induction one obtains the same result
for $[\tau_j,\tau_{j+1}), ~j>j_s$.

\begin{uess}
{\bf [9]} Let $\delta > 0$
be such that $t - h_k (t) < \delta, ~k=1, \dots,
m$.

Then if there exists $K>0$ such that $\mid \! X(t,s) \! \mid < K$
then the equation (1),(2) is stable.

If there exist $N>0,\lambda >0$
such that
$$ \mid \! X(t,s) \! \mid < N e^{- \lambda (t-s)} , $$
then the equation (1),(2) is exponentially stable.
\end{uess}

The proof of the lemma is based on the solution representation
formula (5) and on the fact that $\varphi [h_k (s)] = 0$
for $s > \delta$ .
\vspace{3mm}

Let us study how changing of parameters $A_k,B_j$
influence stability (exponential stability).

Together with the problem (1),(2) we consider
the following one
$$\dot{x} (t) - \sum_{k=1}^m \tilde{A}_k(t) x[h_k(t)]=r(t),
{}~t \in [0, \infty), $$
\begin{eqnarray}
x(\xi)=\varphi(\xi), ~\xi < 0, ~x(\tau_j)=\tilde{B}_j
x(\tau_j -0).    \label{e9}
\end{eqnarray}

\newtheorem{guess}{Theorem}
\begin{guess}
Suppose (a1)-(a4) are satisfied for the problems
(1),(2) and (\ref{e9}) and $0 \leq \tilde{A}_k (t) \leq A_k (t),
{}~0 \leq \tilde{B}_j \leq B_j.$
If the equation (1), (2) is stable (exponentially stable),
then the equation (9) is also stable (exponentially stable).
\end{guess}

{\bf Proof.}
The fundamental function $\tilde{X}(t,s)$ of (9) is a solution
of the problem
$$\dot{x} (t) = \sum_{k=1}^m \tilde{A}_k (t)x[h_k(t)],
{}~t \in [s, \infty), x(s)=1, $$
\begin{eqnarray}
x(\xi)=\varphi(\xi), ~\xi < s, ~x(\tau_j)=\tilde{B}_j
x(\tau_j -0), ~\tau_j > s. \label{e10}
\end{eqnarray}
The problem (\ref{e10}) can be rewritten in the form
$$\dot{x} (t) = \sum_{k=1}^m A_k(t)x[h_k(t)] +
\sum_{k=1}^m \left[ \tilde{A}_k (t) - A_k (t)\right] x[h_k(t)], $$
\begin{eqnarray}
x(s)=1, ~x(\xi)=0, ~\xi <0, \label{e11}
\end{eqnarray}
$$
x(\tau_j)=B_j x(\tau_j -0) + (\tilde{B}_j - B_j)
x(\tau_j -0), ~\tau_j > s.
$$
For solutions of (\ref{e11}) the representation (5), with
the initial point $t=s$, is valid:
\begin{eqnarray}
x(t) = X(t,s) + \sum_{k=1}^m \int_s^t
{X(t,\zeta) [\tilde{A}_k (\zeta) - A_k(\zeta)]
x[h_k(\zeta)] d\zeta} + \nonumber \\
+ \sum_{\tau_j \geq \tau_{j_s}} X(t, \tau_j)
(\tilde{B}_j - B_j) x(\tau_j -0). \label{e12}
\end{eqnarray}

By the corollary of Lemma 2 the solution of (\ref{e10})
is positive, therefore $x(t) >0$.
Besides, $X(t,s) >0,
{}~\tilde{A}_k(\zeta) - A_k(\zeta) \leq 0,
{}~\tilde{B}_j - B_j \leq 0$.
This yields for the solution $x$ of (\ref{e12})
$$x(t) = \tilde{X} (t,s) \leq X(t,s). $$
Referring to Lemma 3 completes the proof.

\section{Main results}

{}~~~~~{\bf 1.}
Consider a special initial value problem
$$
\dot{y} (t) = \sum_{k=1}^m A_k(t)y[h_k(t)],
{}~t \in [\tau_j, \tau_{j+1} ],
$$
\begin{eqnarray}
y(\xi) = 1, \mbox{~if ~~} \xi < \tau_j,
{}~y(\tau_j) = B_j.         \label{e13}
\end{eqnarray}

\begin{guess}
Suppose (a1)-(a4) hold, $A_k(t) \geq 0, ~ B_j \geq 0,$
$$ vrai\sup_{t>0} \int_t^{t+1}\sum_{k=1}^m A_k (s) ds < \infty, $$
there exists $\sigma >0$ such that $\tau_{j+1}
- \tau_j \leq \sigma$ and for any $j$ for the solution $y$
of (\ref{e13}) the inequality
\begin{eqnarray}
y(\tau_{j+1}) \leq 1  \label{e14}
\end{eqnarray}
holds.
Then the equation (1),(2) is stable.
\end{guess}

{\bf Proof.}
In view of Lemma 3 it is sufficient to
demonstrate that the fundamental function $X(t,s)$
is bounded.
Let $s>0$ be fixed and $\tau_{j_s}$ be
the least of all $\tau_j$ such that $\tau_j > s$.
Then $X(t,s)$ as a function of $t$ is a solution
of the problem
$$
\dot{x} (t) = \sum_{k=1}^m A_k(t)x[h_k(t)],
{}~t \in [s, \infty), ~x(s)=1,
$$
\begin{eqnarray}
x(\xi) = 0, \mbox{~if ~~} \xi < s,
{}~x(\tau_j) = B_j x (\tau_j -0), ~\tau_j \geq \tau_{j_s}.
\label{e15}
\end{eqnarray}

Let us prove that it is bounded.
Denote for this solution
$$ \alpha_s = \max_{t \in [s, \tau_{j_s})} x(t). $$
Since $\dot{x} (t) \geq 0 $ in $[s, \tau_{j_s})$
and $x$ is continuous, then $x$
is nondecreasing and $\alpha_s = x(\tau_j -0)$.

Consider the solution of (\ref{e15}) in the interval
$[\tau_{j_s}, \tau_{j_s +1} )$.
For this solution
$$
\dot{x} (t) = \sum_{k=1}^m A_k(t)x(h_k(t)),
{}~t \in [\tau_{j_s}, \tau_{j_s +1} ),
$$
$$
x(\xi) = \varphi(\xi), \mbox{~if ~~} \xi < \tau_{j_s},
{}~x(\tau_{j_s}) = B_{j_s} \alpha_s.
$$
The initial function $\varphi$ is a part of the solution
$x$ up to the point $\tau_{j_s}$.
For this solution  we apply the representation (6)
of the problem (1) without impulses and with
the initial point $t=\tau_{j_s}$
\begin{eqnarray}
x(t) = C(t,\tau_{j_s})x(\tau_{j_s})
+ \sum_{k=1}^m {\int_{\tau_{j_s}}^t {C(t,\zeta)
A_k(\zeta) \varphi [h_k(\zeta)] d\zeta}} ,   \label{e16}
\end{eqnarray}
with $\varphi(h_k(\zeta))=0$, if $h_k(\zeta)>\tau_{j_s}$.

Since $x(\tau_{j_s}) = B_{j_s}\alpha_s$
and $\varphi[h_k(\zeta)] \leq \alpha_s,$
then (\ref{e16}) implies
$$
0 \leq x(t) \leq C(t,\tau_{j_s})B_{j_s} \alpha_s
+ \alpha_s \sum_{k=1}^m {\int_{\tau_{j_s}}^t {C(t,\zeta)
A_k(\zeta) \chi [h_k(\zeta)] d\zeta}} ,
$$

wherein $$\chi(t) = \left\{
\begin{array}{ll}
1,  & \mbox{~if~~} t < \tau_{j_s}, \\
0,  & \mbox{~if~~} t \geq \tau_{j_s}.
\end{array}
\right.
$$

The solution $y$ of (\ref{e13}) for $j=j_s$
can be presented as
$$
y(t) =  C(t,\tau_{j_s})B_{j_s}
+ \sum_{k=1}^m {\int_{\tau_{j_s}}^t {C(t,\zeta)
A_k(\zeta) \chi [h_k(\zeta)] d\zeta}}.
$$
Consequently, for the solution $x$ of the problem (\ref{e15})
the inequality $0 \leq x(t) \leq \alpha_s y(t)$, $t\in
[\tau_{j_s},\tau_{j_s+1}),$ holds.

By (14) $y(\tau_{j_s+1}) \leq 1$.
Thus $x(\tau_{j_s+1}-0) \leq \alpha_s y (\tau_{j_s +1})
\leq \alpha_s$.

The solution of (\ref{e15}) is non-decreasing,
hence $\max_{t \in [\tau_{j_s}, \tau_{j_s +1})}x(t) \leq \alpha_s$.
By induction one obtains that this inequality holds in
any $[\tau_j, \tau_{j+1}), ~j> j_s$, therefore $0<X(t,s) \leq
\alpha_s$.

It remains to prove that $\alpha_s$ is bounded as a function of
$s$.
To this end we consider the problem
$$
\dot{x} (t) = \sum_{k=1}^m A_k(t)x[h_k(t)],
{}~t \in [s, \tau_{j_s}),
$$
\begin{eqnarray}
x(\xi) = 0, \mbox{~if ~~} \xi < s,~x(s)=1. \label{e17}
\end{eqnarray}
If $x$ is a solution of (\ref{e17}), then $\alpha_s =
\max_{t \in [s, \tau_{j_s} )} x(t)$.

In (\ref{e17}) the initial function $\varphi \equiv 0$
and the solution is nondecreasing, thus
$x[h_k(t)] \leq x(t)$.
Therefore for the solution of (\ref{e17}) the inequality
$$ \dot{x} (t) \leq \left( \sum_{k=1}^m A_k (t) \right)
x(t) $$
holds.
The Gronwall-Bellman inequality and $\tau_ {j+1}-\tau_j\leq \sigma$ yield
$$x(t) \leq \exp \left\{ \sum_{k=1}^m
\int_s^t A_k(\zeta) d\zeta \right\}
\leq \exp \left\{ \sum_{k=1}^m
\int_s^{s+\sigma} A_k (\zeta) d \zeta \right\} . $$
Therefore
$$ \alpha_s \leq \exp \left\{ vrai\sup_{s>0} \sum_{k=1}^m
\int_s^{s+\sigma} A_k (\zeta) d \zeta \right\} . $$
Thus $X(t,s)$ can be estimated as
$$ 0 \leq X(t,s) \leq \alpha_s \leq
\exp \left\{ vrai\sup_{s>0} \sum_{k=1}^m
\int_s^{s+\sigma} A_k (\zeta) d \zeta \right\} . $$

It is well known [10] that
$$vrai\sup_{t>0} \int_t^{t+1} |A(s)| ds < \infty$$
implies that for any $\sigma$
$$vrai\sup_{t>0} \int_t^{t+\sigma} |A(s)| ds < \infty,$$
which completes the proof of the theorem.
\vspace{4 mm}

We apply the above theorem for deducing
explicit stability results.
To this end consider (1),(2) for $m=1$ :
$$ \dot{x} (t) = A(t)x[h(t)], ~t \in [0, \infty), $$
\begin{eqnarray}
x(\xi)=\varphi(\xi), ~\xi<0,
{}~x(\tau_j)=B_jx(\tau_j-0). \label{e18}
\end{eqnarray}

The function $h$ is said to satisfy {\bf the separation
condition} in $[\tau_j, \tau_{j+1})$,
if either there exists $t_j \in (\tau_j, \tau_{j+1})$
such that $h(t) \leq \tau_j, ~t \in [\tau_j, t_j),$
and $h(t) \geq \tau_j, ~t \in [\tau_j, \tau_{j+1}),$
or $h(t) \leq \tau_j$ for any $t \in [\tau_j, \tau_{j+1})$.
In the latter case we will assume $t_j = \tau_{j+1}$.
It is to be noted that for $h(t) \equiv t$ any point
$t \in [\tau_j, \tau_{j+1})$ will do.
\begin{guess}
Suppose (a1)-(a4) hold for (\ref{e18}),
$$A(t) \geq 0, ~B_j \geq 0,
{}~vrai\sup_{t>0} \int_t^{t+1} A(s) ds < \infty, $$
there exists $\sigma>0$ such that $\tau_{j+1} - \tau_j
\leq \sigma$ and $h$ satisfies the separation
condition in any interval $[\tau_j, \tau_{j+1})$.

Then the inequality
\begin{eqnarray}
\sup_j \left( B_j +
\int_{\tau_j}^{t_j} A(s) ds \right)
\exp \left\{ \int_{t_j}^{\tau_{j+1}} A(s) ds \right\}
\leq 1,  \label{e19}
\end{eqnarray}
implies the stability of (\ref{e18}).
\end{guess}

{\bf Proof.}
By Theorem 2 one has to estimate solutions of the problem
$$
\dot{y} (t) = A(t)y[h(t)],
{}~t \in [\tau_j, \tau_{j+1}],
$$
\begin{eqnarray}
y(\xi) = 1, \mbox{~if ~~} \xi < \tau_j,~y(\tau_j)=B_j.
\label{e20}
\end{eqnarray}
The solution $y$ of (\ref{e20}) is nondecreasing, therefore
$$ y[h(t)] \leq \left\{
\begin{array}{ll}
1, & \mbox{~if~~} t \in [\tau_j, t_j), \\
y(t), & \mbox{~if~~} t \in [t_j, \tau_{j+1}).
\end{array} \right.
$$
Hence
\begin{eqnarray}
\dot{y}(t) \leq A(t),
\mbox{~if~~} t \in [\tau_j, t_j), \label{e21}
\end{eqnarray}
\begin{eqnarray}
\dot{y}(t) \leq A(t)y(t),
\mbox{~if~~} t \in [t_j, \tau_{j+1}).  \label{e22}
\end{eqnarray}

If (\ref{e21}) holds, then $y(t) \leq B_j + \int_{\tau_j}^t
A(s) ds$. Thus for $t=t_j$
$$y(t_j) \leq B_j + \int_{\tau_j}^{t_j} A(s) ds. $$
Therefore for $t \in [t_j, \tau_{j+1})$ (\ref{e22}) yields
$$y(t) \leq y(t_j) \exp \left\{ \int_{t_j}^t \!\! A(s) ds
\right\}
\leq \left( B_j + \int_{\tau_j}^{t_j} \!\! A(s) ds \right)
\exp \left\{ \int_{t_j}^t \!\! A(s) ds \right\}.$$
Consequently
$$ y(\tau_{j+1}) \leq \left( B_j + \int_{\tau_j}^{t_j}
A(s) ds \right)
\exp \left\{ \int_{t_j}^{\tau_{j+1}} A(s) ds \right\} \leq 1.$$
In view of Theorem 2 the proof is complete.
\vspace{2 mm}

{\bf Remark.}
Since solutions of (1),(2) are bounded in any
finite interval, then it is sufficient for
the hypotheses of Theorems 1 and 2 to be satisfied,
beginning with a certain $j$.
In particular, (\ref{e19}) can be changed by
$$ \lim_{j \rightarrow \infty} \sup_j
\left( B_j + \int_{\tau_j}^{t_j} A(s) ds \right)
\exp \left\{ \int_{t_j}^{\tau_{j+1}} A(s) ds \right\}
\leq 1. $$

Consider the problem
$$ \dot{x}(t) - A(t) x(t-\delta) = r(t), ~t \in [0, \infty) , $$
\begin{eqnarray}
x(\xi) = \varphi (\xi), \mbox{~if~~} \xi < 0,
x(\tau_j) = B_j x(\tau_j -0). \label{e23}
\end{eqnarray}

Denote

$$ \mu_j = \left\{
\begin{array}{lll}
(B_j + \int_{\tau_ j}^{\tau_j +\delta}A(s)ds)\exp\{\int_{\tau_j+\delta}
^{\tau_{j+1}}A(s)ds\}, & \mbox{if} & \delta <\tau _{j+1}-\tau_j , \\
B_j +\int_{\tau_j}^{\tau_{j+1}}A(s)ds, &
\mbox{if} & \delta\geq \tau_{j+1}-\tau_j.
\end{array}
\right.
$$
\vspace{2 mm}

{\bf Corollary}
Suppose
$(a_1), (a_2), (a_4)$  hold for the problem (23),
$$
A(t)\geq 0,~ B_j\geq 0,~ vrai\sup_{t>0} \int_t^{t+1}A(s)ds<\infty,
$$
there exists $\sigma > 0 $ such that $\tau_{j+1}-\tau_j <\sigma $.
Then the inequality $\sup_j \mu_j \leq 1 $ implies the stability of (23).
\vspace{2 mm}

Theorem 3 is easily generalized in a following way.
The function $h$ is said to satisfy {\bf the general separation condition}
in $[\tau_j, \tau_{j+1})$ , if either there exist
$\tau_j=t_j^0 <t_j^1< t_j^2<\dots<  t_j^{k_j} \in[\tau_j, \tau_{j+1})$
such that $ (-1)^{n-1}(\tau_j - h(t))\geq 0,~ t\in [t_j^{n-1}, t_j^n),
n=1,\dots ,k_j, $ or
$h(t)\leq \tau_j$ for any $t\in [\tau_j, \tau_{j+1}).$
In the latter case we will assume $k_j=1, t_j^1=\tau_{j+1}.$

It is to be noted that the general separation condition is not a
restrictive one.
Any function with a finite number of monotonicity
changing on any finite interval satisfies it.

\begin {guess}
Suppose $(a_1)-(a_4)$ hold for (18),
$$
A(t)\geq 0,~ B_j \geq 0,~ vrai\sup_{t>0}\int_t^{t+1}A(s)ds<\infty,
$$
there exists $\sigma >0$ such that $\tau_{j+1}-\tau_j\leq\sigma $
and $h$ satisfies the general separation condition in any interval
$[\tau_j, \tau_{j+1}).$ Then the inequality
$$
\sup_j \left[ \left(\left(B_j+
\!\int_{\tau_j}^{t_j^1}\!\!A(s)ds\right)
\exp\left\{\int_{t_j^1}^{t_j^2}\!\!A(s)ds \right\}+
\right. \right. $$  $$\left. \left.
\!\int_{t_j^2}^{t_j^3}
\!\!\! A(s)ds\right)\exp\left\{\int_{t_j^3}^{t_j^4}\!\!\!A(s)ds
\right\}+\dots
\right]\leq 1
$$
implies the stability of (18).
\end {guess}

Proof is conducted by induction and induction step
from $n$ to $n+1$ is proven as
in Theorem 3.
\vspace{7 mm}

{\bf 2.}
Everywhere above we impose certain conditions
on the delay functions $h_k(t)$.
The following statement gives the stability
result regardless of the delay.
It can be treated as the uniform stabilization
condition for various delay functions.

\begin{guess}
Suppose (a1)-(a4) hold,
$A_k (t) \geq 0, ~k=1, \dots, m,$ and there exists $\sigma
> 0$ such that $\tau_{j+1} - \tau_j \leq \sigma,$
\begin{eqnarray}
q = \sup_k \sup_{t \geq 0} \int_t^{t+ \sigma} A_k (s) ~ds \leq
\frac{1}{m} .  \label{e28}
\end{eqnarray}

If $0 \leq B_j \leq 1-mq, ~j=1,2, \dots, $
then the equation (1),(2) is stable;

if $ 0 \leq B_j \leq 1-mq- \varepsilon,
{}~ \varepsilon > 0, ~j=1,2, \dots,$
and there exists $\rho>0$ such that $h_k (t) \geq t -
\rho$ and $\tau_{j+1}-\tau_j \geq \rho$,
then the equation (1),(2) is exponentially stable.
\end{guess}

\underline{\bf Proof.}
Consider $X(t,s)$ in the segment $[s, \tau_{j_s} )$,
where $\tau_{j_s}$ is the smallest of all $\tau_j > s$.
The function $X(t,s)$ is a solution of the equation (3)
(without impulses).
Since $A_k (t) \geq 0$ then $X(t,s)$ is increasing in $t$.
Thus $X(t,s)$ is majorized by the solution of an ordinary
differential equation
$$ \dot{x} (t) = \left( \sum_{k=1}^m A_k (t) \right) x(t),
{}~x(s) = 1. $$
Therefore by the Gronwall-Bellman inequality
$$X(\tau_{j_s} -0, s) \leq
\exp \left\{ \sum_{k=1}^m \!\int_s^{\tau_{j_s}}\!\! A_k (\zeta) d\zeta
\right\}
\leq $$ $$\exp \left\{ \sum_{k=1}^m \! vrai\sup_{s>0}
\!\int_s^{s+\sigma} \!\!\!\! A_k (\zeta) d\zeta \right\}
\leq e^{mq}
$$
and
$$ X(\tau_{j_s},s) = B_{j_s} X(\tau_{j_s}-0, s) \leq (1-mq)e^{mq}.$$

We shall compare on the interval $[\tau_{j_s}, \tau_{j_s+1} )$
the solutions of two initial value problems
\begin{eqnarray}
\dot{x} (t) = \sum_{k=1}^m A_k (t) x[h_k (t)],
x(\tau_{j_s})=X(\tau_{j_s},s), x(\xi)=X(\xi,s),~\xi < \tau_{j_s},
\label{e29}
\end{eqnarray}
\begin{eqnarray}
\dot{y}(t) = \sum_{k=1}^m A_k (t) e^{mq}, ~y(\tau_{j_s})=X(\tau_{j_s},s).
\label{e30}
\end{eqnarray}

Consider the solution of (\ref{e30}) for $t \in [\tau_{j_s},
\tau_{j_s+1} )$:
$$ y(t) = X(\tau_{j_s},s) + \sum_{k=1}^m \int_{\tau_{j_s}}^t
A_k (\zeta) ~d\zeta ~e^{mq} \leq $$
$$ \leq (1-mq) e^{mq} + mq e^{mq} \leq e^{mq} . $$
Since the solution $x$ of (\ref{e29}) is increasing
for $t \in [\tau_{j_s}, \tau_{j_s+1})$, then
$$ x[h_k(t)] \leq \max \left\{ \sup_{\xi \in [s, \tau_{j_s})}
X(\xi,s), x(t) \right\} \leq \max \left\{ e^{mq}, x(t) \right\} .
$$
Thus the solution $x$ of (\ref{e29}) does not exceed $e^{mg}$.

The solution of (\ref{e29}) in $[\tau_{j_s}, \tau_{j_s+1})$ coincides with
$X(t,s)$.
Then $X(t,s) \leq e^{mq}, ~t \in [s, \tau_{j_s+1} )$.
By induction one easily obtains
$$ X(t,s) \leq e^{mq}, ~ t \in [s, \infty). $$
As $s$ is arbitrary then
$$ \mid \! X(t,s) \! \mid \leq e^{mq}, ~ 0 \leq s \leq t < \infty. $$
By Lemma 3 the equation (1),(2) is stable.

Let  $0 \leq B_j \leq 1-mq - \varepsilon, ~
\varepsilon > 0 $.
Then similarly in the first interval $[s, \tau_{j_s} )$
we have
$$ X(t,s) \leq e^{mq}, ~ t \in [s, \tau_{j_s} ). $$
By the impulsive conditions
$$0 \leq X(\tau_{j_s},s)  \leq B_j e^{mq}
\leq (1-mq-\varepsilon)e^{mq}. $$
By comparing the solutions of (25) and (26)
one obtains
$$X(t,s) \leq (1-mq - \varepsilon)e^{mq} + mqe^{mq} =
(1-\varepsilon) e^{mq}, ~t \in [\tau_{j_s}, \tau_{j_s+1} ). $$

Let us prove by induction that
$$X(t,s) \leq (1-\varepsilon)^i e^{mq}, ~ t \in [\tau_{j_s+i},
\tau_{j_s+i+1} ). $$
We assume
$$X(t,s) \leq (1-\varepsilon)^{i-1} e^{mq}, ~ t \in [\tau_{j_s+i-1},
\tau_{j_s+i} ). $$
In the interval $[\tau_{j_s+i},
\tau_{j_s+i+1} )$
we consider two initial value problems
$$ \dot{x}(t) = \sum_{k=1}^m A_k (t) x[h_k(t)],
{}~x(\tau_{j_s+i})=X(\tau_{j_s+i},s),~
x(\xi)=X(\xi,s),~\xi < \tau_{j_s+i}  $$
and
$$\dot{y(t)} = \sum_{k=1}^m A_k(t)(1 - \varepsilon)^{i-1} e^{mq},
{}~ y(\tau_{j_s+i})= X(\tau_{j_s+i},s). $$
By the hypothesis of the theorem
$\tau_{j_s+i}-\tau_{j_s+i-1} \geq \rho \geq t - h_k(t) > 0$,
consequently in the first equation
$$x[h_k(t)] \leq \max \left\{
\sup_{\xi \in [\tau_{j_s+i-1}, \tau_{j_s+i}) } X(\xi,s), x(t)
\right\}
\leq \max \left\{ (1-\varepsilon)^{i-1} e^{mq}, x(t) \right\}. $$
For the solution $y$ of the second equation
$$y(t) = X(\tau_{j_s+i},s) +
\sum_{k=1}^m \int_{\tau_{j_s+i}}^t A_k(\xi)
(1-\varepsilon)^{i-1} e^{mq} d\xi  \leq $$
$$ (1-mq -\varepsilon)(1-\varepsilon)^{i-1}e^{mq}+
mq(1-\varepsilon)^{i-1}e^{mq}=$$
$$
(1-\varepsilon)^i e^{mq} < (1-\varepsilon)^{i-1} e^{mq}. $$

{}From here the right hand side of the first equation
does not exceed the right hand side of the second equation.
Since the initial values are equal,
then
$$X(t,s)=x(t) \leq y(t) \leq
(1-\varepsilon)^i e^{mq}, ~ t \in [\tau_{j_s+i},
\tau_{j_s+i+1} ), $$
which completes the induction step.

By the hypothesis $\tau_{j+1}-\tau_j \geq \rho$,
therefore
$$ X(t,s) \leq (1-\varepsilon)^{(t-s)/  \rho -1} e^{mq}, $$
and $X(t,s)$ has an exponential estimate
\begin{eqnarray}
|X(t,s) | \leq N e^{-\lambda (t-s)}, ~\mbox{with}
{}~N=e^{mq}/(1-\varepsilon), ~\lambda = -\frac{1}{\rho}
\ln (1-\varepsilon).
\label{e31}
\end{eqnarray}

Then by Lemma 3 the equation (1),(2) is exponentially stable.
The proof of the theorem is complete.

\vspace{3 mm}

{\bf 3.}
Now we proceed to the stabilization of (1),(2) without assuming
$A_k (t) \geq 0$.
Let write these coefficients in the form
$$A_k(t) = A_k^+ (t) - A_k^- (t), $$
where $a^+ = \max \{ a,0 \}, ~ a^- = a^+ - a $.

Consider an auxiliary equation
\begin{eqnarray}
\dot{x} (t) - \sum_{k=1}^m H_k (t) x[g_k(t)] = f(t), ~t \geq 0,
x(\xi) = \varphi (\xi), \mbox{~~if~} \xi < 0 ; \nonumber \\
x(\tau_j) = B_jx(\tau_j-0),~j=1,\dots.
\label{e32}
\end{eqnarray}

Suppose that for the equation (\ref{e32}) the hypotheses of Lemma 1
are satisfied.
Then for this problem the solution representation formula (5)
holds with a certain fundamental function $R(t,s)$.

The operator
$$ (Rz)(t) = \int_0^t R(t,s) z(s) ds $$
is said to be the Cauchy operator of the problem (\ref{e32}).

By ${\bf L}_{\infty}$ we denote a space of functions
$x: [0, \infty) \rightarrow {\bf R}$ Lebesgue measurable and
essentially bounded on $[0, \infty)$,
with the norm
$$\parallel x \parallel_{{\bf L}_{\infty}} =
vrai\sup_{t>0} \mid \! x(t) \! \mid . $$

In sequel the following proposition will be necessary;
it is a corollary from the results of [5].

\begin{uess}
{\bf [5]}
Suppose (a1)-(a4) hold, $ \sup_{j>0} \mid \! B_j \! \mid < \infty$,
there exist $\rho > 0,\delta > 0$ such that
$\tau_{j+1} - \tau_j \geq \rho, ~t-h_k(t)<\delta$
and the fundamental function $R(t,s)$ of the problem (\ref{e32})
has an exponential estimate
$$ \mid \! R(t,s) \! \mid < N e^{- \lambda (t-s) }, $$
with $N>0,~\lambda > 0$.
If the operator ${\cal L} R : {\bf L}_{\infty} \rightarrow {\bf
L}_{\infty}$ is invertible,
where the operator ${\cal L}$ is defined by the left-hand side
of (3), then the equation (1),(2) is exponentially stable.
\end{uess}

\begin{guess}
Suppose (a1)-(a4) hold, $B_j \geq 0, ~ \sup_j B_j < \infty$,
there exist $\rho > 0, ~\delta > 0$ and $\sigma >0$ such that
$\rho \leq \tau_{j+1} - \tau_j \leq \sigma, ~ t - h_k(t) <
\delta$,
$$
q=\sup_k \sup_{t \geq 0} \int_t^{t+\sigma} A_k^+ (s) ds <
\infty, ~~   B_j \leq 1 - mq - \varepsilon,
{}~\varepsilon > 0,
$$
\begin{eqnarray}
vrai\sup_{t \geq 0} \left( \sum_{k=1}^m A_k^- (t) \right)
< \frac{\lambda}{N},   \label{e33}
\end{eqnarray}
where $N$ and $\lambda$ are defined in (\ref{e31}).

Then the equation (1),(2) is exponentially stable.
\end{guess}

{\bf Proof.}
We prove this theorem by applying Lemma 4.
Let $R$ be the Cauchy operator of the equation
$$\dot{x}(t) - \sum_{k=1}^m A_k^+(t) x(h_k(t)) = f(t), ~t \geq 0,
$$
$$
x(\xi) = \varphi(\xi), ~ \mbox{~~if~} \xi < 0, $$
with impulsive conditions (2), $R(t,s)$
be a fundamental function of this equation.

The proof of Theorem 4 yields that the fundamental function
satisfies the estimate
\begin{eqnarray}
0 \leq R(t,s) \leq N e^{- \lambda (t-s)}, \label{e34}
\end{eqnarray}
where $N$ and $\lambda$ are defined by (\ref{e31}).
We write the operator ${\cal L}$ of the problem (1),(2)
in the form
$$ ({\cal L} x )(t) = \dot{x} (t) -
\sum_{k=1}^m A_k^+ (t) x[h_k(t)] +
\sum_{k=1}^m A_k^- (t) x[h_k (t)], ~ t \geq 0, $$
$$x(\xi) = 0, ~ \xi < 0 . $$

For the operator ${\cal L} R$ we have
${\cal L} R = E + T$, where $E$ is the identity operator,
$$ (Tz)(t) = \sum_{k=1}^m A_k^- (t)
\int_0^{h_k^+ (t)} R(h_k (t), s) z(s) ~ds. $$
Thus the estimate (\ref{e34}) yields
$$ \parallel T \parallel_{{\bf L}_{\infty} \rightarrow
{\bf L}_{\infty}} \leq
\frac{N}{\lambda} vrai\sup_{t>0}
\sum_{k=1}^m A_k^- (t) . $$

Therefore (\ref{e33}) gives $\parallel T \parallel < 1 $.
Consequently the operator ${\cal L}R : {\bf L}_{\infty}
\rightarrow {\bf L}_{\infty}$ is invertible,
so all the hypotheses of Lemma 4 are satisfied,
which completes the proof of the theorem.

\end{document}